\documentclass[draft,11pt,english]{article}

\usepackage{amsfonts,amsmath,amsthm,amscd,amssymb,latexsym,verbatim,
caption2}

\theoremstyle{plain}

\newtheorem{lemma}{Lemma}[section]

\theoremstyle{definition}

\newcommand{\keywords}{\textbf{Key words. }\medskip}
\newcommand{\subjclass}{\textbf{MSC 2010. }\medskip}
\renewcommand{\abstract}{\textbf{Abstract. }\medskip}

\numberwithin{equation}{section}
\def\MYim{\mathop{\rm Im}\limits}
\def\MYre{\mathop{\rm Re}\limits}

\begin{document}

\title {Darlington's Theorem for Rational Nevanlinna Functions in Several Variables}

\author{M. F. Bessmertny\u{\i}}




\maketitle

\begin{abstract}
Every rational Nevanlinna function in $n$ variables is a Cayley  inner function in $n + 1$ variables with one variable fixed in the upper half-plane. Darlington's theorem reduces the problem of synthesis of a system with losses to the problem of synthesis without losses.
\end{abstract}
\medskip




    \section*{Introduction}\label{s:0}
    Let $\Pi^{d}=\{z\in\mathbb{C}^{d}\mid \MYim{z_{1}}>0,\ldots,\MYim{z_{d}}>0\}$ be an open upper poly-half-plane.
    \emph{The Nevanlinna class}
    $\mathcal{N}_{d}^{m\times m}$ consists of $\mathbb{C}^{m\times m}$-valued functions $f(z_{1},\ldots,z_{d})$ holomorphic in $\Pi^{d}$ and satisfying the condition
      \begin{equation}\label{eq0.1}
        \MYim f(z)=(f(z)-f(z)^{\ast})/2i\geq0\quad \mbox{for} \quad z\in\Pi^{d},
      \end{equation}
  where $^{\ast}$ means transition to the Hermitian conjugate matrix.
  We will consider only \emph{rational function} of the Nevanlinna class, as well as a subclass $\mathcal{IN}_{d}^{m\times m}$ of Nevanlinna rational functions that are \emph{Cayley inner}. For such functions
      \begin{equation}\label{eq0.2}
      (f(x)-f(x)^{\ast})/2i=0\quad \mbox{for} \quad x\in\mathbb{R}^{d}.
      \end{equation}
  The latter means that such a function is the double Cayley transform (over the variables and over the matrix values) of an inner
  function on the unit polydisk  $\mathbb{D}^{d}=\{z\in\mathbb{C}^{d}\mid |z_{1}|<1,\ldots,|z_{d}|<1\}$.
  The structure of scalar rational inner functions in polydisk is considered in \cite{uj6}.

  We can consider the class $\mathcal{R}_{d}^{m\times m}$ of positive real ($\mathbb{R}^{m\times m}$-valued) functions $z(\lambda)$ for which
     $$
     \MYre z(\lambda)\geq 0 \quad\text{for}\quad \MYre \lambda_{1}>0,\ldots,\MYre \lambda_{d}>0
     $$
  and, accordingly, the subclass $\mathcal{IR}_{d}^{m\times m}$ of Cayley inner positive real functions.

  The impedances of passive electrical circuits are one-variable rational positive real functions.  Various methods are known for realizing a such function as the impedance of an electric circuit. In particular, the well-known Darlington theorem (\cite{uj03}, see also \cite{uj5}) admits the following equivalent mathematical formulation.
  \medskip

  \emph{For any scalar one-variable rational positive real function $z(\lambda)$  there exists a rational positive real Cayley inner second-order matrix function}
     \begin{equation}\label{eq0.3}
     Z(\lambda)=
       \begin{pmatrix}
       a(\lambda) & b(\lambda) \\
       c(\lambda) & d(\lambda)
       \end{pmatrix}
     \end{equation}
  \emph{such that
     \begin{equation}\label{eq0.4}
     z(\lambda)=a(\lambda)-b(\lambda)\left[\,d(\lambda)+r\,\right]^{-1}c(\lambda),
     \end{equation}
  where $r$ is a non-negative constant}.
  \medskip

  \noindent
  It is easy to see that
     \begin{equation}\label{eq0.5}
     g(\lambda,\mu)=a(\lambda)-b(\lambda)\left[\,d(\lambda)+\mu\,\right]^{-1}c(\lambda)
     \end{equation}
  is the positive real Cayley inner function in the right bi-half-plain $\Omega^{2}=\{\lambda,\mu\in\mathbb{C}^{d}\mid \MYre \lambda>0, \MYre \mu>0\}$ and
     \begin{equation}\label{eq0.6}
     z(\lambda)=g(\lambda,r),\quad\text{where}\quad r>0.
     \end{equation}

  We generalize this result to rational functions in several variables of the Nevanlinna class.
  \medskip

  \noindent
  \textbf{Main Theorem} (\textbf{Generalized Darlington's theorem}). \emph{For any rational $\mathbb{C}^{m\times m}$-valued function $f(z)\in\mathcal{N}_{d}^{m\times m}$ there exists a rational Cayley inner function  $g(z,z_{d+1})\in\mathcal{IN}_{d+1}^{m\times m}$ such that}
    \begin{equation}\label{eq0.7}
    f(z)=g(z,i).
    \end{equation}

  \noindent
  This theorem allows one to obtain a representation of an arbitrary rational function of the Nevanlinna class from a representation of Cayley inner functions of the Nevanlinna class (see \cite{uj01} -- \cite{uj2}).

    \section{Preliminaries and Proof of the Main Theorem}\label{s:1}

    Let $f(z)=\{f_{ij}(z)\}_{i,j=1}^{m}$ be a rational $\mathbb{C}^{m\times m}$-valued function and $q(z)$ be the common denominator of the functions $f_{ij}(z)$. Then $f_{ij}(z)=p_{ij}(z)/q(z)$, where $p_{ij}(z)$ are polynomials. The function $f(z)$ will be written in the form
    $f (z)=P(z)/q(z)$, where $P(z)=\{p_{ij}(z)\}_{i,j=1}^{m}$ is a $\mathbb{C}^{m\times m}$-valued polynomial and $q(z)$ is a scalar polynomial. In fact, division $P(z)/q(z)$ is the standard operation of multiplying of the matrix $P(z)$ by the number $q(z)^{-1}$.

    We say a polynomial $p(z)$ with complex coefficients is \emph{stable}\footnote{This terminology may differ from the designations of other authors.} if it is nonzero for any $z\in \Pi^{d}$. A stable polynomial with real coefficients will be called \emph{a real stable}. The ring of the polynomials with real coefficients is denoted by $\mathbb{R}[z_{1},\ldots,z_{d}]$.
    \medskip

    To prove Main Theorem, we need two lemmas.

    \begin{lemma}\label{lem1.1} \emph{(\cite{uj4}, Corollary 5.5)}.
    Let $p(z)+iq(z)\neq0$ where $p(z),q(z)\in\mathbb{R}[z_{1},\ldots,z_{d}]$, and let $z_{d+1}$ be a new indeterminate. Then following are equivalent.

    \emph{(i)} $p(z)+iq(z)$ is stabile,

    \emph{(ii)} $p(z)+z_{d+1}q(z)$ is a real stabile polynomial in $d+1$ variables,

    \emph{(iii)} all nonzero polynomials in the pencil
        $$
        \{\alpha p(z)+\beta q(z)\mid \alpha,\beta\in\mathbb{R}\}
        $$
    are real stable.
    \end{lemma}

    \begin{lemma}\label{lem1.2}\emph{(\cite{uj3}, Lemma 2.2 and \cite{uj4}, Theorem 5.4)}.
    Let $p(z),q(z)\in\mathbb{R}[z_{1},\ldots,z_{d}]$ be coprime and $p(z)\not\equiv 0$, $q(z)\not\equiv 0$. Then all non-zero polynomials in the pencil
      $$
      \{\alpha p(z)+\beta q(z)\mid \alpha,\beta\in\mathbb{R}\}
      $$
    are real stable if and only if
      $$
      \MYim p(z)/q(z)\geq 0\quad\text{or}\quad \MYim p(z)/q(z)\leq 0
      $$
    whenever $\MYim z_{k}>0$ for all $1\leq k\leq d$.
    \end{lemma}
    \medskip

    \noindent
    \emph{Proof of Main Theorem}.
    Let $f(z)=P(z)/q(z)\in \mathcal{N}_{d}^{m\times m}$ be a rational  function and let $P(z)$, $q(z)$ be a coprime polynomials. Then
      $$
      f_{\eta}(z)=\frac{\eta P(z)\eta^{\ast}}{q(z)}
      $$
    is the scalar function of the Nevanlinna class for every row vector $\eta\in \mathbb{C}^{m}$.
    The function $f_{\eta}(z)+\alpha$ also belongs to $\mathcal{N}_{d}^{1\times 1}$ for any $\alpha\in\mathbb{R}$.
    Since $P(z)$ and $q(z)$ are coprime polynomials, then the numerator $\alpha q(z)+\eta P(z)\eta^{\ast}$ of $f_{\eta}(z)+\alpha$ is a stable polynomial with complex coefficients.

    We set
     \begin{equation}\label{eq1.1}
     P_{1}(z)=[P(z)-P(\overline{z})^{\ast}]/2i,\quad P_{2}(z)=[P(z)+P(\overline{z})^{\ast}]/2,
     \end{equation}
     \begin{equation}\label{eq1.2}
     q_{1}(z)=[q(z)-\overline{q(\overline{z})}]/2i,\,\quad
     q_{2}(z)=[q(z)+\overline{q(\overline{z})}]/2.\quad
     \end{equation}
    Clearly, the $\mathbb{C}^{m\times m}$-valued polynomials $P_{1}(z),P_{2}(z)$ have Hermitian coefficients. Then the polynomial $\eta P_{1}(z)\eta^{\ast}, \eta P_{2}(z)\eta^{\ast}$ and $q_{1}(z), q_{2}(z)$ have a real coefficients.

    Using (\ref{eq1.1}), (\ref{eq1.2}), we obtain
      \begin{multline}\nonumber
      \alpha q(z)+\eta P(z)\eta^{\ast}=\alpha (iq_{1}(z)+q_{2}(z))+\eta (iP_{1}(z)+P_{2}(z))\eta^{\ast}=\\
      i(\alpha q_{1}(z)+\eta P_{1}(z)\eta^{\ast})+
      (\alpha q_{2}(z)+\eta P_{2}(z)\eta^{\ast}).\qquad
      \end{multline}
    Let $z_{d+1}$ be a new indeterminate. By Lemma \ref{lem1.1}, for every $\eta\in\mathbb{C}^{m}$
      \begin{multline}\nonumber
      z_{d+1}(\alpha q_{1}(z)+\eta P_{1}(z)\eta^{\ast})+
      (\alpha q_{2}(z)+\eta P_{2}(z)\eta^{\ast})=\\
      \alpha(z_{d+1}q_{1}(z)+q_{2}(z))+
      \eta (z_{d+1}P_{1}(z)+P_{2}(z))\eta^{\ast}
      \end{multline}
    is a real stable polynomial in $d+1$ variables. Since $\eta$ is an arbitrary vector, all nonzero polynomials of the pencil
      $$
      \alpha(z_{d+1}q_{1}(z)+q_{2}(z))+
      \beta(\eta (z_{d+1}P_{1}(z)+P_{2}(z))\eta^{\ast}),\quad \alpha,\beta\in\mathbb{R}
      $$
    are real stable. By condition
       $$
       \MYim \left.\frac{\eta(z_{d+1}P_{1}(z)+P_{2}(z))\eta^{\ast}}
       {z_{d+1}q_{1}(z)+q_{2}(z)}\right|_{z_{d+1}=i}=
       \MYim \eta f(z)\eta^{\ast}\geq 0
       \quad\text{for all}\quad z\in\Pi^{d}.
       $$
    Therefore, by Lemma \ref{lem1.2},
      $$
      \MYim \frac{\eta(z_{d+1}P_{1}(z)+P_{2}(z))\eta^{\ast}}{z_{d+1}q_{1}(z)+q_{2}(z)}\geq 0
      \quad\text{for}\quad (z,z_{d+1})\in \Pi^{d+1}.
      $$
    Since $\eta\in\mathbb{C}^{m}$ is an arbitrary vector and
      $$
      \MYim\frac{x_{d+1}P_{1}(x)+P_{2}(x)}{x_{d+1}q_{1}(x)+q_{2}(x)}=0
      \quad\text{for}\quad (x,x_{d+1})\in \mathbb{R}^{d+1},
      $$
    then
      $$
      g(z,z_{n+1})=\frac{z_{d+1}P_{1}(z)+P_{2}(z))}{z_{d+1}q_{1}(z)+q_{2}(z)}
      \in\mathcal{IN}_{d+1}^{m\times m}.
      $$
    Moreover, from (\ref{eq1.1}), (\ref{eq1.2}) follows $g(z,i)=f(z)$, what was required.
    \qed

   \end{document}